\documentclass[sts,authoryear,preprint]{imsart}

\RequirePackage[OT1]{fontenc}
\usepackage{amsthm,amsmath,amsfonts,natbib}
\RequirePackage[colorlinks,citecolor=blue,urlcolor=blue]{hyperref}

\startlocaldefs
\numberwithin{equation}{section}
\theoremstyle{plain}

\endlocaldefs

\begin{document}

\begin{frontmatter}
\title{Is profile likelihood a true likelihood? An argument in favor}
\runtitle{Profile likelihood}

\begin{aug}
\author{\fnms{Oliver J.} \snm{Maclaren}\ead[label=e1]{oliver.maclaren@auckland.ac.nz}}

\runauthor{Maclaren}

\affiliation{Department of Engineering Science, University of Auckland, Auckland, New Zealand}

\address{Department of Engineering Science, The University of Auckland, Auckland 1142, New Zealand \printead{e1}.}

\end{aug}

\begin{abstract}
Profile likelihood is the key tool for dealing with nuisance parameters in likelihood theory. It is often asserted, however, that profile likelihood is not a `true' likelihood. One implication is that likelihood theory lacks the generality of e.g. Bayesian inference, wherein marginalization is the universal tool for dealing with nuisance parameters. Here we argue that profile likelihood has as much claim to being a true likelihood as a marginal probability has to being a true probability distribution. The crucial point we argue is that a likelihood function is naturally interpreted as a maxitive possibility measure: given this, the associated theory of integration with respect to maxitive measures delivers profile likelihood as the direct analogue of marginal probability in additive measure theory. Thus, given a background likelihood function, we argue that profiling over the likelihood function is as natural (or as unnatural, as the case may be) as marginalizing over a background probability measure. The connections to Bayesian inference can also be further clarified with the introduction of a suitable logarithmic distance function, in which case the present theory can be naturally described as `Tropical Bayes' in the sense of tropical algebra.
\end{abstract}


\end{frontmatter}

\section{Introduction}
Consider the opening sentence from the entry on profile likelihood in the Encyclopedia of Biostatistics \citep{Aitkin2005-re}:

\begin{quote}The profile likelihood is not a likelihood, but a likelihood maximized over nuisance parameters given the values of the parameters of interest.\end{quote}

Numerous similar assertions that profile likelihood is not a `true' likelihood may be found throughout the literature and various textbooks, and is apparently the accepted viewpoint of the statistical community. Importantly, this includes the `pure' likelihood literature, which generally accepts a lack of systematic methods for dealing with nuisance parameters, while still recommending profile likelihood as the most general, albeit `ad-hoc', solution \citep[see e.g.][]{Royall1997-ai,Rohde2014-tz,Edwards1992-nj,Pawitan2001-xm}. Similarly, recent monographs on characterizing statistical evidence presents favorable opinions of the likelihood approach but criticize the lack of general methods for dealing with nuisance parameters \citep{Aitkin2010-ea,Evans2015-kg}. The various justifications given, however, appear to the present author to rather vague and unconvincing. For example, suppose we modified the above quotation to refer to marginal probability instead of profile likelihood:

\begin{quote}A marginal probability is not a probability, but a probability distribution integrated over nuisance variables given the values of the variables of interest.\end{quote}

The above would be a perfectly fine characterization of a marginal probability if the ``not a probability, but'' part was dropped, i.e. 

\begin{quote}A marginal probability is a probability distribution integrated over nuisance variables given the values of the variables of interest.\end{quote}

Simply put: the fact that a marginal probability is obtained by integrating over a `background' probability distribution does not prevent the marginal probability from being a true probability. The crucial observation in the case of marginal probability is that \textit{integration over variables takes probability distributions to probability distributions}. 

The purpose of the present article is to point out that \textit{there is an appropriate notion of integration over variables that takes likelihood functions to likelihood functions via maximization}. This notion of integration is based on the idea of \textit{idempotent analysis}, wherein one replaces a standard algebraic operation such as addition in a given mathematical theory with another basic algebraic operation, defining a form of `idempotent addition', to obtain a new analogous, self-consistent theory \citep{Maslov1992-wo,Kolokoltsov1997-rc}. In this case one simply replaces the usual `addition' operations, including the usual (Lebesgue) integration, with `maximization' operations, including taking supremums, to obtain a new, `idempotent probability theory'. Maximization in this context is understood algebraically as an idempotent addition operation, hence the terminology. While perhaps somewhat exotic at first sight, this idea finds direct applications in e.g. large deviation theory \citep{Puhalskii2001-qk} and, most relevantly, possibility theory, fuzzy set theory and pure-likelihood-based decision theory \citep{Dubois1997-js,Cattaneo2013-ub,Cattaneo2017-sd}. A popular special instance of idempotent mathematics is so-called `tropical mathematics' in which multiplication is also converted to a new algebraic operation, here addition \citep[see e.g.][]{Speyer2009-ic,Akian1996-oj,Litvinov2007-li,Pachter2004-sp,Bernhard2000-lr}. That is, the basic `addition' and `multiplication' operations in tropical algebra are interpreted as $(\text{max},+)$, respectively, instead of the usual $(+,\times)$. With the introduction of a logarithmic distance in likelihood theory, multiplication of likelihoods becomes addition of log-likelihoods and we are naturally led to a `Tropical Bayesian' interpretation of (log) profile likelihoods. This provides a formal foundation for the usual intuitive interpretation of (negative) log-likelihoods as `cost' measures.

The present argument is not, of course, without objections. In particular, acceptance or rejection of the present interpretation depends on what one believes the key properties of likelihood should be; this is, perhaps surprisingly, not without significant controversy \citep{Bayarri1992-zd,Bjornstad1996-ju,Bayarri1988-uq}. Thus we end with a discussion of various potential objections, including a discussion of some properties one might want a general notion of `likelihood' to satisfy and whether the present interpretation does or does not satisfy these. Despite potential conflicts with some frequentist, evidential and/or Bayesian considerations, we believe that the present interpretation is a clear, self-consistent and suitable foundational concept for `pure' likelihood theory \citep[particularly that developed by][]{Edwards1992-nj}, and/or for what we propose to call `Tropical Bayes'.


\section{Likelihood as a possibility measure}
\label{sec:meth}
Though apparently not well known in the statistical literature, likelihood theory is known in the wider literature on uncertainty quantification to have a natural correspondence to possibility theory rather than to probability theory \citep{Dubois1997-js,Cattaneo2013-ub,Cattaneo2017-sd}. This has perhaps been obscured by the usefulness of likelihood methods as tools in probabilistic statistical inference. It is not our intention to review this wider literature in detail here \citep[see e.g.][for more]{Dubois1997-js,Cattaneo2013-ub,Cattaneo2017-sd,Augustin2014-su,Halpern2017-mc}, but to simply point out the implications of this correspondence. In particular, likelihood theory interpreted as a possibilistic, rather than probabilistic theory can be summarized as: 

\begin{quote}Probability theory with addition replaced by maximization.\end{quote}

As indicated above, this is sometimes known as, for example, `idempotent measure theory', `maxitive’ measure theory or `possibility' theory, among other names \citep[see e.g.][for more]{Dubois1997-js,Cattaneo2013-ub,Cattaneo2017-sd,Augustin2014-su,Halpern2017-mc,Maslov1992-wo,Kolokoltsov1997-rc,Puhalskii2001-qk}. This correspondence perhaps explains the preponderance of maximization methods in likelihood theory, including the methods of maximum likelihood and profile likelihood. 

The most important consequence of this perspective is that the usual Lebesgue integration with respect to an additive measure, as in probability theory, becomes, in likelihood/possibility theory, a different type of integration, defined with respect to a maxitive measure. Again, the key point is simply that addition operations (including summation and integration) are replaced by maximization operations (or taking supremums in general).

For completeness, we contrast the key axioms of possibility theory with those of probability theory. Given a set of possibilities of $\Omega$, assumed to be discrete for the moment for simplicity, and for two discrete sets of possibilities $A,B \subseteq \Omega$ the key axioms of elementary possibility theory are \citep{Halpern2017-mc}:

\begin{equation}
\begin{aligned}
&\text{poss}(\emptyset) = 0\\
&\text{poss}(\Omega) = 1\\
&\text{poss}(A \cup B) = \text{max}\{\text{poss}(A),\text{poss}(B)\}
\end{aligned}
\end{equation}

which can be contrasted with those of elementary probability theory:

\begin{equation}
\begin{aligned}
&\text{prob}(\emptyset) = 0\\
&\text{prob}(\Omega) = 1\\
&\text{prob}(A \cup B) = \text{sum}\{\text{prob}(A),\text{prob}(B)\}
\end{aligned}
\end{equation}

where $A$ and $B$ are required to be disjoint in the probabilistic case, but this is not strictly required in the possibilistic case.

Given a `background' or `starting' likelihood measure, likelihood theory can be developed as a self-contained theory of possibility, where derived distributions are manipulated according to the first set of axioms above. This is entirely analogous to developing probability theory from a background measure, with derived distributions manipulated according to the second set of axioms. As our intention is to consider methods for obtaining derived distributions by `eliminating' nuisance parameters, we need not consider here where the starting measure comes from (but see the Discussion).

To make the correspondences of interest clear in what follows, we first present probabilistic marginalization as a special case of a pushforward measure or, equivalently, as a special case of a general (not necessarily 1-1) change of variables. We then consider the possibilistic analogues.

\section{Pushforward probability measures and the delta function method for general changes of variable}
\label{sec:verify}

Given a probability measure $\mu$ over a random variable $x \in \mathbb{R}^n$ with associated density $\rho$, define the new random variable $t=T(x)$ where $T:\mathbb{R}^n \rightarrow \mathbb{R}^m$. This variable is distributed according to the pushforward measure $T\star \mu$, i.e. $t \sim T\star \mu$. 

The density of $t$, here denoted by $q = T\star\rho$, is conveniently calculated via the delta function method which is valid for arbitrary changes of variables (not necessarily 1-1):

\begin{equation}
q(t) = [T\star \rho](t) = \int \delta(t-T(x))\rho(x)dx.
\end{equation}

As a side point, we note that this method of carrying out arbitrary transformations of variables is standard in statistical physics \citep[see e.g.][]{Van_Kampen1992-ik}, but is apparently less common in statistics \citep[see the articles][aimed at highlighting this method to the statistical community]{Au1999-oh,Khuri2004-mo}.

\subsection{Marginalization via the delta function method}
The above means that we can interpret marginalization to a component $x_1$, say, as a special case of a (non-1-1) deterministic change of variables via:

\begin{equation}
\rho(x_1) = \int \delta(x_1-\text{proj}_{X_1}(x))\rho(x)dx,
\end{equation}

where $\text{proj}_{X_1}(x)$ is simply the projection of $x$ to its first coordinate. Thus marginalization can be thought of as the pushforward under the projection operator and as a special case of a general (not necessarily 1-1) change of variables $t = T(x)$.

\section{Profile likelihood as marginal possibility and an extension to general changes of variable}
As we have repeatedly stressed above, likelihood theory interpreted as a possibilistic, and hence maxitive, measure theory simply means that ‘addition’ operations such as the usual Lebesgue integration are replaced by ‘maximization’ operations such as taking the supremum.

Consider first then the analogue of a marginal probability density, which we will call a marginal possibility distribution and denote by $L_p$. Starting from a `background' likelihood measure $L(x)$ we `marginalize' in the analogous manner to before:

\begin{equation}
L_p(x_1) = \text{sup}\{\delta(x_1-\text{proj}_{X_1}(x))L(x)\} = 
\text{sup}_{\{x | \text{proj}_{X_1}(x) = x_1\}}\{L(x)\}.
\end{equation}

This is again simply the pushforward under the projection operator, but here under a different type of `integration' - i.e. the operation of taking a supremum. Of course, this is just the usual profile likelihood for $x_1$. 



As above, we need not be restricted to ‘marginal’ possibility distributions: we can consider arbitrary functions of the parameter $t = T(x)$. This leads to an analogous pushforward operation of $L(x)$ to $L_p(t)$ that we denote by $\star_p$:

\begin{equation}
L_p(t) = [T\star_p L](t) = \text{sup}\{\delta(t-T(x))L(x)\} = \text{sup}_{\{x | T(x) = t\}}\{L(x)\}
\end{equation}

which again corresponds to the usual definition of profile likelihood.

\section{A simple example comparing marginal probability and marginal possibility}
Here we consider a simple example illustrating the difference between probabilistic and possibilistic reasoning, in particular under marginalization/non-1-1 changes of variable. 

Suppose you have three suspects in a crime. Through some means or another you decide on the following `plausibility' distribution, where plausibility is used here as a general umbrella term for probability and/or possibility reasoning: suspect one has plausibility 0.4, while the other two suspects each have plausibility 0.3. You also know that suspect one was wearing a red hat at the time of the crime while the other two were wearing blue hats. 

According to the above, under a \textit{probabilistic} interpretation, the most probable \textit{perpetrator} is suspect one (who wore a red hat); but the most probable \textit{hat color} of the perpetrator is blue (with probability 0.3 + 0.3 = 0.6). This is a consequence of the \textit{additivity} of probability theory and the non-1-1 change of variables in going from suspects to hat colors.

On the other hand, if you interpret the given plausibility numbers as a \textit{possibility} distribution, then according to standard possibility theory the most possible suspect is suspect one and the most possible hat colour is now red, i.e. is the hat color of the most possible suspect, suspect one. Similarly, this is a consequence of the \textit{maxitivity} of possibility theory.

The difference can be made more extreme given a large number of `other' suspects, each with low plausibility but sharing some common property that the main suspect lacks. Again, these results are a simple consequence of how additivity and maxitivity, respectively, interact with non-1-1 changes of variable (here: person to hat color). 

We believe that there are reasonable situations where additivity is desirable, \textit{but also} reasonable situations in which maxitivity might be preferred. This is a subject worth further debate. We note, however, that a \textit{relative probability} approach to the problem of statistical evidence, such as that presented in \cite{Evans2015-kg} comes to similar conclusions to that of a possibility approach (Michael Evans, personal communication).

\section{Strength of evidence, distances and `Tropical Bayes'}
As noted in \cite{Evans2015-kg}, it is perhaps less controversial to hold that likelihood gives a \textit{qualitatively} reasonable relative \textit{ordering} of preference for parameter values in light of data than it is to hold \citep[e.g.][]{Royall1997-ai} that it provides a \textit{quantitative} measure of relative support.

To make some progress towards addressing this distinction, we consider how to define a suitable notion of \textit{distance} that respects - but is distinct from - a given \textit{qualitative ordering}. Notions of statistical distance are common in the statistical literature \citep[see e.g.][and references therein]{Basu2011-wa}; here, however, we follow the ideas developed by \cite{Tarantola2006-cb} of \textit{quality spaces} and distances defined in these. This leads naturally to the idea of pure likelihood theory as a form of what we propose to call `Tropical Bayes', where the meaning of this term is discussed below.

In particular, given the \textit{ordering} induced by a likelihood function (and/or profile likelihood function):

\begin{equation}\theta_1 \text{ is preferred to } \theta_2 \text{ iff } L(\theta_1) > L(\theta_2), 
\end{equation}

we can define a \textit{likelihood distance} via

\begin{equation}
D_{L}(\theta_1,\theta_2) = \vert \log \frac{L(\theta_2)}{L(\theta_1)}\vert = \vert \log \frac{L(\theta_1)}{L(\theta_2)}\vert.
\end{equation}

This distance has the properties of being symmetric, additive and zero iff $L(\theta_1) = L(\theta_2)$. \cite{Tarantola2006-cb} argues that this notion of distance is widely applicable for many types of qualitative orderings. In the present case it is, of course, just the well-known log-likelihood ratio function. We propose then that, accepting that the likelihood gives a natural qualitative preference or plausibility \textit{ordering}, the log-likelihood then gives a natural \textit{distance} in this `qualitative space'. There remains, however, a choice of logarithm base and/or a choice of arbitrary distance scale factor; thus we can't fully remove some of the `qualitative' features associated with pure likelihood theory without a further choice of reference. One natural choice might be to take the minimum distance to a fully saturated model, i.e. one which can fit the data perfectly, in which case one would be interested in how much `fit' to trade-off against parsimony considerations \citep{Edwards1992-nj}.

Interestingly, the combination of replacing addition operations by maximization and then working in log-space (wherein multiplication becomes addition) corresponds to completing the `tropicalization' of probability theory: moving from an algebraic structure in terms of $(+,\times)$ to one in terms of $(\text{max},+)$. This is the subject of `tropical algebra', which also goes by the name `max-plus' algebra, and is a popular special instance of idempotent mathematics with applications to decision theory, uncertainty quantification, statistical inference and optimization \citep[see e.g.][for some relevant starting points in this area]{Speyer2009-ic,Akian1996-oj,Litvinov2007-li,Pachter2004-sp,Bernhard2000-lr}. A natural interpretation of negative log-likelihood functions in this context is as `cost measures'; these have also been termed `Maslov measures', due to their origins in Maslov's idempotent probability theory \citep{Akian1996-oj,Bernhard2000-lr}. These analogies are explored in detail by \cite{Akian1996-oj}, where the natural analogue of a random variable is a decision variable, the analogue of a Markov chain is a Bellman chain (i.e. the Bellman equation from the subject of dynamic programming) and so on.

Finally, however, we note that even if profile likelihood is accepted as the natural analogue of marginal probability, the \textit{evidential} interpretation of profile likelihood may still have difficulties; this is discussed further below.

\section{Discussion}
\subsection{Objections to profile likelihood}
As discussed, it is frequently asserted that profile likelihood is not a true likelihood \citep{Aitkin2005-re,Royall1997-ai,Pawitan2001-xm,Rohde2014-tz,Evans2015-kg}. Common reasons include: that it is obtained from a likelihood via maximization \citep{Aitkin2005-re}, that it is not based directly on observable quantities \citep{Royall1997-ai,Pawitan2001-xm,Rohde2014-tz} and that it lacks particular repeated sampling properties \citep{Royall1997-ai,Cox1994-ww}. 

None of the above objections appear to the present author to apply to the following: given a starting or `background' likelihood function, profile likelihood satisfies the axioms of possibility theory, in which the basic additivity axiom of probability theory is replaced by a maxitivity axiom. Profile likelihood is simply the natural possibilistic counterpart to marginal probability, where additive integration is replaced by a maxitive analogue. We thus argue that, if marginal probability is a `true' probability, then profile likelihood should likewise be considered a `true' likelihood, at least when likelihood theory is interpreted in a possibilistic manner. Negative log-likelihood functions can then be naturally interpreted as cost measures in the sense of tropical mathematics.

\subsection{Fixed data}
Regarding the second two objections mentioned above: observable quantities and repeated sampling properties, it is important to note that the given data must be held fixed to give a consistent background likelihood over which to profile. Given fixed data one has a fixed possibility measure and thus can consider `marginal' - i.e. profile - likelihoods. In contrast, repeated sampling will produce a distribution of such possibility measures, and these may or may not have good frequentist properties. None of this is in contrast to marginal probability: changing the distribution over which we marginalize changes the resulting marginal probability. Of course, despite this caveat, profile likelihood often does have good repeated sampling properties \citep{Royall1997-ai,Cox1994-ww} and also plays a key role in frequentist theory, though we do not discuss this further here. One consequence is that our conception of profile likelihood does not generally satisfy properties such as zero expectation of the associated score function \citep{Cox1994-ww,Pawitan2001-xm}. These are, however, properties dependent on particular repeated sampling notions such as `unbiasedness', and hence more properly considered as frequentist concepts. The present approach is more suitable for those seeking a non-probabilistic `plausibility' measure, as induced by data that are considered fixed once observed. 

\subsection{Why?}
A natural question, perhaps, is why worry about whether profile likelihood is a true likelihood? One answer is that profile likelihood is a widely used tool but is often dismissed as `ad-hoc' or lacking proper justification. This gives the impression that, for example, likelihood theory is lacking in comparison with e.g. Bayesian theory in terms of systematic methods for dealing with nuisance parameters. By understanding that profile likelihood does in fact have a systematic basis in terms of possibility theory practitioners and students can better understand and reason about a widely popular and useful tool. Understanding the connection to \textit{possibilistic} as opposed to probabilistic reasoning may also help explain why profile likelihood has emerged as a particularly promising method of \textit{identifiability analysis} \citep{Raue2009-sl}, where identifiability is traditionally a prerequisite for probabilistic analysis. Of course, as indicated, the price of accepting profile likelihood as a `true' likelihood is an interpretation in terms of pure likelihood theory, and this makes the connections to repeated sampling properties more complicated. We see no need however, to restrict oneself to one perspective on statistical inference - the present possibilistic view can \textit{complement} other approaches such as frequentist statistics or Bayesian statistics. Furthermore, this analogy opens strong connections between likelihood theory and the optimization literature; the foundations of such connections have already been explored by e.g. \cite{Akian1996-oj,Bernhard2000-lr} and provide a natural link to pure likelihood decision theory as developed by \cite{Cattaneo2013-ub}.

\subsection{Ignorance}
The possibilistic interpretation of likelihood also helps understand the representation of ignorance. While probabilistic ignorance is not preserved under \textit{arbitrary} changes of variables (e.g. non-1-1 transformations), even in the discrete case, possibilistic ignorance is in the following sense: if we take the maximum likelihood over a set of possibilities, such as $\{x\ |\ T(x) = t\}$ for each $t$, rather than summing them, a flat `prior likelihood' \citep{Edwards1969-vc,Edwards1992-nj} over $x$ becomes a flat prior likelihood over $t$. On the other hand, a flat prior probability over $x$ in general becomes non-flat over $t$ under non-1-1 changes of variable. Thus a profile prior likelihood has what, in many cases, may be desirable properties as a representation of prior ignorance \cite[see the discussion in][for more on likelihood and the representation of ignorance]{Edwards1969-vc,Edwards1992-nj}. This difference in transformation properties was also illustrated in our simple example comparing the probabilistic and possibilistic analysis of criminal evidence. As noted there, however, the \textit{relative} probabilistic approach \textit{a la} \cite{Evans2015-kg}, reaches conclusions closer to the possibilistic analysis, compared to the conclusions of the `absolute' probabilistic analysis (Michael Evans, personal communication).


\subsection{Point function or set function?}
Likelihood is traditionally considered a point function as opposed to a set function; this is also related to controversy over defining likelihood functions for so-called composite hypotheses \citep[see e.g.][]{Edwards1992-nj,Royall1997-ai}. Authors such as \cite{Basu2012-wm} have argued, \textit{contra} e.g. Fisher, that likelihood could be directly extended to a set function. \cite{Basu2012-wm} further developed the argument that this set function could be taken as additive - we are more inclined, here at least, to consider the first possibility, and reject the second. A number of other authors have also considered the question of composite hypotheses, in particular in the context of defining evidence \citep[see e.g.][]{Zhang2013-cu,Blume2013-nr,Bickel2012-nb}.

We have attempted to avoid the issue of set functions/composite hypotheses somewhat by instead using the concept of a non-1-1 transformation of variables. This allows us to consider the likelihood of subsets of the full/background parameter space based on an indexing statistic, i.e. by using subsets defined via $\{x\ |\ T(x) = t\}$. This approach is based on what amounts to equality constraints, leaving out subsets defined via inequality constraints. It may be desirable to further relax this and simply consider likelihood directly as a set function defined via

\begin{equation}
L_p(A) = \text{sup}_{x \in A}\{L(x)\}
\end{equation}

for $A \subseteq X$. This allows for inequality constraints such as those in $A = \{x\ |\ T(x) \leq t\}$.

We leave consideration of this approach to future work. Presumably, however, one could recover the present approach by considering some notion of minimal and/or extremal sets of equality constraints, e.g. by restricting attention to those inequality constraints that are \textit{active} during the profiling/maximization procedure, and hence those that are reduced to binding equality constraints. The interpretation of negative log-likelihoods as cost measures may also be helpful here.

\subsection{Evidence}
One of the key issues to consider when deciding whether to accept profile likelihoods as `true' likelihoods is whether they can play the same role that `full' likelihoods play in defining \textit{evidential} measures \citep{Royall1997-ai,Aitkin2010-ea,Evans2015-kg,Zhang2013-cu,Blume2013-nr,Bickel2012-nb}. Mathematically, it appears clear that profile likelihood is entirely analogous to marginal probability; it is less clear whether - or under what circumstances - one should use marginal (whether maxitive or additive) measures in defining \textit{evidence}. We believe that this applies equally to the Bayesian approach. A way forward from here would be to separate the questions: first accept profile likelihood as a `marginal' possibility measure, and then investigate under what circumstances marginal measures can be given further evidential interpretations. We suspect that the answer may require additional concepts and/or assumptions like those used in the causal inference literature to separate spurious marginal associations from `true' causation \citep{Pearl2009-jp,Pearl2009-qh}. That is, we suspect that `evidence' may be better defined in causal terms than in either purely probabilistic or purely possibilistic terms. As such, the question of whether or not profile likelihood is a `true' likelihood should be independent of whether it plays the role of an evidential measure, unless the definition of likelihood is itself explicitly supplemented with causal assumptions.

\section{Conclusions}
\label{sec:conc}
We have argued that profile likelihood has as much claim to being a true likelihood as a marginal probability has to being a true probability distribution. In the case of marginal probability, integration over variables takes probability distributions to probability distributions, while in the case of likelihood, \textit{maximization takes likelihood functions to likelihood functions}. Maximization can be considered in this context as an alternative (idempotent) notion of integration, and a likelihood function as a maxitive possibility measure. There are some conflicts with both Bayesian and frequentist considerations, however: lack of additivity and lack of some repeated sampling properties, respectively. In our view, these conflicts are not necessarily an issue, as neither additivity nor repeated sampling properties such as unbiasedness are beyond objections. Instead we argue that the present approach gives a self-consistent theory suitable for \textit{possibilistic} statistical analysis, with a well-defined method of treating nuisance parameters, and which continues in the tradition of `pure' likelihood theories. The connection of profile likelihoods to evidential interpretations appears subtle (as is, we believe, the connection of marginal probabilities to evidence); our view is that this issue should be explored further in the context of formulating additional \textit{causal} properties that an evidence measure should satisfy, such as those required to classify marginal correlations into `spurious' and `true' causal relationships. Finally, taking profile likelihood seriously as a `true' likelihood leads naturally to the idea of `Tropical Bayesian Inference', a subject yet to be properly explored by the statistical community.

\section*{Acknowledgements}
The author would like to thank Michael Evans, Marco Cattaneo, Yudi Pawitan, Alexandre Patriota, Christian Robert and Anthony Edwards for useful comments and/or discussions.

\bibliographystyle{imsart-nameyear}

\bibliography{Bibliography-MM-MC}

\end{document}